\documentclass[a4paper, 12pt]{amsart}

\usepackage{amsmath, amsfonts, amssymb, amsthm, graphics,graphicx}

\newtheorem{thm}{Theorem}[section]
\newtheorem{lemma}[thm]{Lemma}
\newtheorem{cor}[thm]{Corollary}

\providecommand{\aut}{\mathop{\rm Aut \,}\nolimits}

\providecommand{\gmod}{\mathop{\rm mod}\nolimits}

\renewcommand{\\}{\vspace{3mm}}

\begin{document}

\title[\bf Orbital digraphs of primitive groups]{\bf Orbital digraphs of infinite primitive permutation groups}

\author{\bf Simon M. Smith}

\thanks{Address: 59 Netherwood Road, London W14 0BP. Formally, Mathematical Institute, University of Oxford, Oxford.}
\thanks{Email: simon.smith@chch.oxon.org}

\subjclass{20 B 15; 05 C 25}

\date{\today}

\begin{abstract}
If $G$ is a group acting on a set $\Omega$ and $\alpha, \beta \in
\Omega$, the digraph whose vertex set is $\Omega$ and whose
arc set is the orbit $(\alpha, \beta)^G$ is called an {\em
orbital digraph} of $G$. Each orbit of the stabiliser $G_\alpha$
acting on $\Omega$ is called a {\it suborbit} of $G$.

A digraph is {\em locally finite} if each vertex is adjacent to at
most finitely many other vertices. A locally finite digraph $\Gamma$
has more than one end if there exists a finite set of vertices $X$
such that the induced digraph $\Gamma \setminus X$ contains at least
two infinite connected components; if there exists such a set
containing precisely one element, then $\Gamma$ has {\em
connectivity one}.

In this paper we show that if $G$ is a primitive permutation group
whose suborbits are all finite, possessing an orbital digraph with
more than one end, then $G$ has a primitive connectivity-one
orbital digraph, and this digraph is essentially unique. Such digraphs
resemble trees in many respects, and have been fully characterised
in a previous paper by the author.

\end{abstract}

\maketitle

\newpage
\section{Preliminaries}

Throughout this paper, $G$ will be a group of permutations of an
infinite set $\Omega$. A transitive
group $G$ is {\it primitive} on $\Omega$ if the only
$G$-invariant equivalence relations admitted by $\Omega$ are the trivial and universal
relations. A transitive group $G$ is
said to act {\em regularly} on $\Omega$ if $G_\alpha = 1$ for each
$\alpha \in \Omega$.

If $\alpha \in \Omega$ and $g \in G$, we denote the image of
$\alpha$ under $g$ by $\alpha^g$. Following this notation, all
permutations will act on the right. The set of images of $\alpha$
under all elements of $G$ is called an {\it orbit} of $G$, and is
denoted by $\alpha^G$. There is a natural action of $G$ on the
$n$-element subsets and $n$-tuples of $\Omega$ via $\{\alpha_1,
\ldots, \alpha_n\}^g := \{\alpha_1^g, \ldots, \alpha_n^g\}$ and
$(\alpha_1, \ldots, \alpha_n)^g := (\alpha_1^g, \ldots,
\alpha_n^g)$ respectively.

If $\alpha \in \Omega$ we denote the stabiliser of $\alpha$ in $G$
by $G_\alpha$, and if $\Sigma \subseteq \Omega$ we denote the
setwise and pointwise stabilisers of $\Sigma$ in $G$ by
$G_{\{\Sigma\}}$ and $G_{(\Sigma)}$ respectively.\\

A {\it digraph} $\Gamma$ will be a directed graph without multiple edges or
loops; it is a pair $(V\Gamma,
A\Gamma)$, where $V\Gamma$ is the set of {\it vertices} and
$A\Gamma$ the set of {\it arcs} of $\Gamma$. The set $A\Gamma$
consists of ordered pairs of distinct elements of $V\Gamma$. Two
vertices $\alpha, \beta \in V\Gamma$ are {\it adjacent} if either
$(\alpha, \beta)$ or $(\beta, \alpha)$ lies in $A\Gamma$. A digraph
is {\em locally finite} if every vertex is adjacent to at most
finitely many vertices. All paths will be undirected, unless
otherwise stated. A digraph is {\em infinite} if its vertex set is
infinite. The {\em distance} between two
connected vertices $\alpha$ and $\beta$ in $\Gamma$ is denoted by
$d_\Gamma(\alpha, \beta)$.

A {\it graph} is a pair $(V\Gamma,
E\Gamma)$, in which $E\Gamma$ is a set of unordered pairs of elements
of $V\Gamma$, each of which is known as an {\it edge} of $\Gamma$.
The results presented in this paper concern digraphs,
and are thus applicable to graphs.

The {\em automorphism group of $\Gamma$} will be denoted by $\aut
\Gamma$. A graph or digraph is {\em primitive} if $\aut \Gamma$ is primitive
on the set $V\Gamma$, and is {\em automorphism-regular} if $\aut
\Gamma$ acts regularly on $V\Gamma$.\\

If $G$ is a group of permutations of a set $\Omega$ and $\alpha,
\beta \in \Omega$, then the digraph with vertex set
$\Omega$ and arc set $(\alpha, \beta)^G$ is called an {\em
orbital digraph of $G$}. Note such digraphs are necessarily
arc-transitive. The orbital digraph $(\Omega, (\alpha, \beta)^G)$
is {\em non-diagonal} if $\alpha$ and $\beta$ are distinct.
Naturally, the graph $(\Omega, \{\alpha, \beta\}^G)$ is called an {\it orbital graph} of $G$.

The connected components of any digraph $\Gamma$ are
$\aut \Gamma$-invariant equivalence classes on the vertex set of $\Gamma$;
thus any non-diagonal primitive orbital digraph is connected.

\section{Ends of graphs}

A concept that is fundamental to the arguments contained herein is
that of an {\em end}. A {\em half-line} $L$ of a graph or digraph $\Gamma$ is
a one-way infinite cycle-free path $x_1 x_2 x_3 \ldots$ in
$\Gamma$. The vertex $x_1$ is commonly known as the {\em root} of
$L$. There is a natural equivalence relation on the set of
half-lines of $\Gamma$: two half-lines $L_1$ and $L_2$ are
equivalent if there exist an infinite number of pairwise-disjoint
paths connecting a vertex in $L_1$ to a vertex in $L_2$. The
equivalence classes of this relation are called the {\em ends of
$\Gamma$}. An end is called {\em thick} if it contains infinitely
many disjoint half-lines, otherwise it is called {\em thin}. By
convention, a finite graph or digraph has no ends.

\begin{thm} {\normalfont(\cite[Corollary 1.5]{watkins:survey})}
\label{thm:no_of_ends_of_graphs} If $\Gamma$
is an infinite locally finite connected vertex-transitive graph
then $\Gamma$ has $1$, $2$ or $2^{\aleph_0}$ ends. \qed \end{thm}

Given a digraph $\Gamma$, it is possible to construct an associated graph
$\Gamma^\prime$, whose vertex set is $V\Gamma$, with two vertices being
adjacent in $\Gamma^\prime$ if and only if they are adjacent in $\Gamma$.
Clearly, the ends of a digraph and its associated graph are the same. Thus, the above
theorem is also true of digraphs.

Much of what is currently known about graphs with more than one
end is underpinned by the work of Dicks and Dunwoody in their
seminal text \cite{dicks&dunwoody:gaog}. An accessible summary of
their results pertaining to locally finite graphs can be found in
\cite{moller:survey}.

The {\it connectivity} of an infinite connected graph or digraph $\Gamma$ is
the smallest possible cardinality of a subset $W$ of $V\Gamma$ for
which the induced digraph $\Gamma \setminus W$ is not connected. A
{\it lobe} of $\Gamma$ is a connected subgraph that is maximal
subject to the condition it has connectivity strictly greater than
one. If $\Gamma$ has connectivity one, then the vertices $\alpha$
for which $\Gamma \setminus \{\alpha\}$ is not connected are
called the {\it cut vertices} of $\Gamma$.

Let $V_1$ be the set of cut
vertices of $\Gamma$, and let $V_2$ be a set in
bijective correspondence with the set of lobes of $\Gamma$. Define $T$ to
be a bipartite digraph whose parts are $V_1$ and $V_2$. Two
vertices $\alpha \in V_1$ and $x \in V_2$ are adjacent in $T$ if
and only if $\alpha$ is contained in the lobe of $\Gamma$
corresponding to $x$; in this case there is a pair of arcs
between them, one in each direction.
This construction yields a tree, which is called the {\it
block-cut-vertex tree} of $\Gamma$. If $\Gamma$ has
connectivity one and block-cut-vertex tree $T$, then for any group $G$
acting on $\Gamma$, there is a natural action on $T$.\\

Two digraphs $\Gamma_1$ and $\Gamma_2$ are {\it quasi-isometric} if
there exists a map $\phi: V\Gamma_1 \rightarrow V\Gamma_2$ and
there exist constants $a, b > 0$ such that for all $\alpha, \beta
\in V\Gamma$ we have
\[\frac{1}{a} \, d_{\Gamma_1}(\alpha, \beta) - b \leq d_{\Gamma_2}(\phi(\alpha), \phi(\beta)) \leq a.
d_{\Gamma_1}(\alpha, \beta) + b.\] In \cite[Proposition
1]{moller:ends_of_graphs_2}, M\"{o}ller shows if two locally
finite connected graphs are quasi-isometric, then there is a
natural bijection between their ends. This is also true of digraphs.

\begin{lemma} \label{lemma:orbital_graphs_have_same_ends} If $G$ is
a primitive group of permutations of an infinite set
$\Omega$ and $\Gamma_1$ and $\Gamma_2$ are two non-diagonal locally finite
orbital digraphs of $G$ then there is a natural bijection between
the ends of $\Gamma_1$ and the ends of $\Gamma_2$. \end{lemma}

\begin{proof} Suppose both 
$\Gamma_1 = (\Omega, (\alpha_1, \beta_1)^G)$ and $\Gamma_2 = (\Omega, (\alpha_2, \beta_2)^G)$ 
are locally finite orbital digraphs of $G$.
Put $m_1:=d_{\Gamma_1}(\alpha_2, \beta_2)$ and $m_2:=d_{\Gamma_2}(\alpha_1, \beta_1)$,
and let $a:= \mathrm{max \,} \{m_1, m_2\}$. Since both $\Gamma_1$ and $\Gamma_2$ are connected,
any arc in $\Gamma_2$ corresponds to a path of length $m_1$ in $\Gamma_1$, and any arc in
$\Gamma_1$ corresponds to a path of length $m_2$ in $\Gamma_2$. Thus, if $\alpha, \beta \in \Omega$
then $d_{\Gamma_1}(\alpha, \beta) \leq m_1.d_{\Gamma_2}(\alpha, \beta) \leq a.d_{\Gamma_2}(\alpha, \beta)$,
and $d_{\Gamma_2}(\alpha, \beta) \leq m_2.d_{\Gamma_1}(\alpha, \beta) \leq a.d_{\Gamma_1}(\alpha, \beta)$.
Therefore $(1 / a) d_{\Gamma_1}(\alpha, \beta)
\leq d_{\Gamma_2}(\alpha, \beta) \leq a.d_{\Gamma_1}(\alpha, \beta)$,
so $\Gamma_1$ and $\Gamma_2$ are quasi-isometric.\end{proof}

Henceforth, no distinction will be made between an end of $\Gamma_1$ and its corresponding end in $\Gamma_2$. Furthermore, we define a {\em permutation-end} of $G$ to be an end of an orbital digraph of $G$. By the above lemma, the set of permutation-ends of $G$ is equal to the set of ends of a non-diagonal orbital digraph of $G$, and is independent of the orbital digraph chosen.\\

\section{The canonical orbital digraph}

In \cite{moller:primitivity}, M\"{o}ller shows every connected
locally finite primitive graph with more than one end
is closely related to one with connectivity one.

\begin{thm} {\normalfont(\cite[Theorem 15]{moller:survey})} \label{thm:jung_watkins_moller}
If $\Gamma$ is a connected locally finite primitive
graph with more than one end then there exist vertices $\alpha,
\beta \in V\Gamma$ such that the graph $(V\Gamma,
\{\alpha, \beta\}^{\aut \Gamma})$ has connectivity one and each
lobe has at most one end. \qed \end{thm}

Since any infinite vertex-transitive connectivity-one digraph has
infinitely many ends, the following may be deduced.

\begin{cor} \label{corollary:no_2_ended_graphs} If $\Gamma$ is a
locally finite primitive arc-transitive digraph then $\Gamma$ has
$0$, $1$ or $2^{\aleph_0}$ ends. \qed \end{cor}

If $G$ is a group of permutations of a set $\Omega$, then the orbits of
the action of the stabiliser $G_\alpha$ on $\Omega$
are called the {\em $\alpha$-suborbits} of $G$.
The cardinality of a suborbit is called a {\em subdegree} of
$G$.

Recall that the set of permutation-ends of a
primitive group with no infinite subdegree is
equal to the set of ends of any one of its non-diagonal
orbital digraphs. A consequence of the above corollary therefore, is
that infinite primitive groups with no infinite subdegree either have
precisely one permutation-end, or they have uncountably many.

We can use Theorem~\ref{thm:jung_watkins_moller} to tell us a
great deal about the orbital digraphs of infinite primitive
permutation groups.

\begin{thm} \label{thm:small_connectivity1} If $G$ is a primitive group
of permutations of an infinite set $\Omega$, with more than one
permutation-end and no infinite subdegree, then $G$ has a locally
finite orbital digraph with connectivity one in which each lobe has
at most one end.
\end{thm}

The proof of this theorem relies on the following five lemmas. Let
$G$ be a primitive group of permutations of an infinite set
$\Omega$.

\begin{lemma} \label{lemma:G_has_connectivity_one_orbital_graph} If $G$ has
a locally finite orbital digraph with more than one end then $G$ has
an orbital digraph with connectivity one.\end{lemma}

\begin{proof} Suppose $G$ has a locally finite orbital digraph $\Gamma$ with more than one end.
Let $\Gamma^\prime$ be the graph associated with $\Gamma$.
By Theorem~\ref{thm:jung_watkins_moller},
we may find vertices $\alpha, \beta \in \Omega$ such that the
graph $(\Omega, \{\alpha, \beta\}^{\aut \Gamma^\prime})$ has
connectivity one. Now $G$ is a primitive subgroup of $\aut
\Gamma$, and thus $\aut \Gamma^\prime$, so the graph $(\Omega,
\{\alpha, \beta\}^G)$ is a connected subgraph of $(\Omega,
\{\alpha, \beta\}^{\aut \Gamma^\prime})$, and therefore has connectivity
one. Hence, the digraph $\Gamma:= (\Omega, (\alpha,
\beta)^G)$ also has connectivity one. \end{proof}

\begin{lemma} {\normalfont(\cite[Lemma 2.1.]{me:InfPrimDigraphs})} \label{lemma:block_stabiliser_is_primitive}
If $G$ acts primitively on the vertices, and transitively on the arcs, of a connectivity-one digraph
$\Gamma$, and $\Lambda$ is any lobe of $\Gamma$ with at least three vertices, then the subgroup of $\aut \Lambda$ induced by the setwise stabiliser $G_{\{\Lambda\}}$ acts
primitively on the vertices of $\Lambda$. \qed \end{lemma}

The following three lemmas are generalisations of several
observations made in \cite{moller:primitivity}. Although
M\"{o}ller restricts his attention to automorphism groups of
graphs, for our purposes his arguments require only
trivial modification.

\begin{lemma} \label{lemma:same_block}
If $\Gamma$ is a locally finite connectivity-one orbital digraph of
$G$ and $\alpha, \beta \in V\Gamma$ lie in the same lobe
$\Lambda$ of $\Gamma$, then a lobe of $\Lambda^\prime :=
(V\Lambda, (\alpha, \beta)^{G_{\{\Lambda\}}})$ is a lobe of
$\Gamma^\prime:= (V\Gamma, (\alpha, \beta)^G)$. \qed \end{lemma}

Let $\Gamma$ be a locally finite orbital digraph of $G$ with more
than one end. Fix $\alpha, \beta \in V\Gamma$ such that
$\Gamma^{(1)}:= (V\Gamma, (\alpha, \beta)^G)$ has connectivity
one. Let $\Lambda_1$ be a lobe of $\Gamma^{(1)}$ containing
$\alpha$, and let $G^{(1)}:= G_{\{\Lambda_1\}}$. This group acts
primitively on $V\Lambda_1$ by
Lemma~\ref{lemma:block_stabiliser_is_primitive}. If $\Lambda_1$
has more than one end then find $\beta_2 \in \Lambda_1$ such that
$\Gamma^{(2)}:= (V\Lambda_1, (\alpha, \beta_2)^{G^{(1)}})$ has
connectivity one. Let $\Lambda_2$ be a lobe of $\Gamma^{(2)}$
containing $\alpha$ and put $G^{(2)} := G_{\{\Lambda_2\}}$. Again
we note this group acts primitively on $V\Lambda_2$ by
Lemma~\ref{lemma:block_stabiliser_is_primitive}. For $i \geq 2$,
if $\Lambda_i$ has more than one end, find $\beta_{i+1} \in
\Lambda_i$ such that $\Gamma^{(i+1)}:= (V\Lambda_i, (\alpha,
\beta_{i+1})^{G^{(i)}})$ has connectivity one. Let $\Lambda_{i+1}$
be a lobe of $\Gamma^{(i+1)}$ containing $\alpha$ and put
$G^{(i+1)}:= G_{\{\Lambda_{i+1}\}}$. By
Lemma~\ref{lemma:block_stabiliser_is_primitive}, this group acts
primitively on $V\Lambda_{i+1}$.

\begin{lemma} \label{lemma:value_of_n} There is an $n \geq 1$
such that $\Lambda_n$ has at most one end. \qed \end{lemma}

\begin{lemma} \label{lemma:lambda_i_block_of_gamma_prime} For all integers
$i \geq 1$, the
digraph $\Lambda_i$ is a lobe of the orbital digraph
$(V\Gamma, (\alpha, \beta_i)^G)$. \qed \end{lemma}

\begin{proof} [Proof of Theorem~\ref{thm:small_connectivity1}] Suppose $G$ has
more than one permutation-end and no infinite subdegree.

If $\Gamma$ is a non-diagonal orbital digraph of $G$, then $\Gamma$
has more than one end. Furthermore, if $(\alpha, \beta) \in
A\Gamma$, and $(\gamma, \alpha) \in A\Gamma$, then $\alpha$ is
adjacent to at most $|\beta^{G_\alpha}| + |\gamma^{G_\alpha}|$
vertices in $\Gamma$. Since $G$ has no infinite subdegrees, this
sum is finite, and since $G$ acts transitively on $V\Gamma$, this
digraph, and indeed all orbital digraphs of $G$, are locally finite.

Thus, by Lemma~\ref{lemma:G_has_connectivity_one_orbital_graph},
$G$ has a locally finite connectivity-one orbital digraph $\Gamma$.
Using the notation described above, if we set
$\Gamma^{(1)}:=\Gamma$ then, by Lemma~\ref{lemma:value_of_n},
there exists a positive integer $n$ such that $\Lambda_n$ has at
most one end. Put $\Gamma':=(V\Gamma, (\alpha, \beta_n)^G)$. By
Lemma~\ref{lemma:lambda_i_block_of_gamma_prime}, $\Lambda_n$ is a
lobe of $\Gamma^\prime$, and since $\Lambda_n \not = \Gamma'$, we
see $\Gamma'$ has connectivity one. Furthermore, $G$ acts
arc-transitively on $\Gamma'$, so every lobe of $\Gamma'$ is
isomorphic to $\Lambda_n$, and thus has at most one end.
\end{proof}

The structure of primitive directed digraphs with connectivity one
is well known.

\begin{thm} {\normalfont(\cite[Theorem 2.6.]{me:InfPrimDigraphs})}
\label{theorem:block_stabiliser_not_regular} Let $G$ be a
vertex-transitive group of automorphisms of a connectivity-one
digraph $\Gamma$ whose lobes have at least three vertices. If $G$
acts primitively on $V\Gamma$ and $\Lambda$ is some lobe of
$\Gamma$, then $G_{\{\Lambda\}}$ is primitive and not regular on
$V\Lambda$. \qed \end{thm}

Recall that a digraph $\Gamma$ is automorphism-regular if $\aut
\Gamma$ acts regularly on $V\Gamma$.

\begin{thm} {\normalfont(\cite[Theorem 3.3.]{me:InfPrimDigraphs})} \label{thm:directed_connectivity_one_graphs}
If $\Gamma$ is a
vertex-transitive digraph with connectivity one, then it is
primitive if and only if the lobes of $\Gamma$ are primitive but
not automorphism-regular, pairwise isomorphic and each has at
least three vertices. \qed \end{thm}

Combining Theorem~\ref{thm:small_connectivity1},
Theorem~\ref{theorem:block_stabiliser_not_regular} and
Theorem~\ref{thm:directed_connectivity_one_graphs} we obtain the
following.

\begin{thm} \label{thm:connectivity1}
If $G$ is a primitive group of permutations of an infinite set
$\Omega$ with more than one permutation-end and no infinite
subdegree, then $G$ has a locally finite orbital digraph $\Gamma$
with connectivity one, whose lobes are primitive but not
automorphism-regular, are pairwise isomorphic, have at least three
vertices and at most one end. Furthermore, if $\Lambda$ is a lobe
of $\Gamma$, then $G_{\{\Lambda\}}$ acts primitively but not
regularly on $V\Lambda$. \qed \end{thm}

If $G$ is a primitive group of permutations of an infinite set
$\Omega$ and $\Gamma$ is a locally finite orbital digraph with
connectivity one whose lobes have at most one end, then we shall
call $\Gamma$ a {\em canonical orbital digraph of $G$}. If $\Gamma$
is such a digraph, and $\Lambda$ is a lobe of $\Gamma$ containing
precisely one end $\epsilon$, then $\epsilon$ will be called a
{\em lobe-end} of $\Gamma$. If an end $\epsilon$ is not a
lobe-end of $\Gamma$, it will be called a {\em tree-end}, as it
corresponds to an end of the block-cut-vertex tree of $\Gamma$.

Observe that $G$ permutes the lobe-ends of $\Gamma$ transitively.
Indeed, $G$ acts arc-transitively on $\Gamma$, and each arc
belongs to a unique lobe of $\Gamma$, so $G$ permutes the lobes,
and therefore the lobe-ends, transitively.

Since $\Gamma$ has connectivity one, any thick ends of $\Gamma$
must be lobe-ends. Thus, if $\Gamma$ has a thick end then $G$
acts transitively on the thick ends of $\Gamma$.

\section{Uniqueness}

Recall that if $G$ is a primitive group and $\Gamma_1$ and
$\Gamma_2$ are locally finite orbital digraphs of $G$ then there is a natural
bijection between the ends of $\Gamma_1$ and the ends of
$\Gamma_2$.

A {\em geodesic} between two connected vertices of a graph or digraph $\Gamma$ is a
path of minimal length between them. If $T$ is a tree, and
$\alpha, \beta \in VT$, there is exactly one geodesic in $T$
between $\alpha$ and $\beta$, which we denote by $[\alpha,
\beta]_T$. If $\epsilon$ is an end of $T$, then there is precisely
one half-line in $T$ that lies in $\epsilon$ with initial vertex
$\alpha$. We denote this half-line by $[\alpha, \epsilon)_T$. A
connected subgraph $C$ of $\Gamma$ is said to {\it contain
the end $\epsilon$} if, for every half-line $(\alpha_i)_{i \in
\mathbb{N}}$ lying in $\epsilon$, we have $\alpha_i \in  VC$ for
all sufficiently large $i$.

\begin{lemma} \label{lemma:thick_ends_go_to_thick_ends}
Suppose $\Gamma_1$ and $\Gamma_2$ are canonical orbital digraphs of
a primitive group $G$. If $\epsilon$ is a lobe-end of $\Gamma_1$
then $\epsilon$ is a lobe-end of $\Gamma_2$. \end{lemma}

\begin{proof} Let $T_1$ and $T_2$ be the lobe-cut-vertex trees of $\Gamma_1$
and $\Gamma_2$ respectively. Suppose, for a contradiction, that
there is a lobe-end $\epsilon$ of $\Gamma_1$ that is a tree-end
of $\Gamma_2$. Since $\epsilon$ is a tree-end of $\Gamma_2$, it is
an end of $T_2$. As it is a lobe-end of $\Gamma_1$, there exists
a lobe $\Lambda_1$ of $\Gamma_1$ whose end is $\epsilon$, and a
vertex $x \in VT_1 \setminus V\Gamma_1$ corresponding to the lobe
$\Lambda_1$.

\begin{figure}
\centering
\input{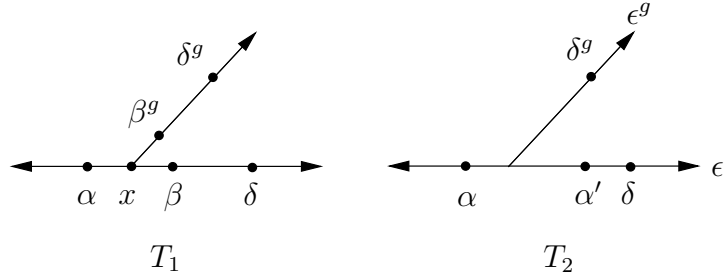}
\caption{The trees $T_1$ and $T_2$}
\label{T1_and_T2}
\end{figure}

The following argument is illustrated in Figure~\ref{T1_and_T2}.
Fix $\alpha \in V\Lambda_1$, and let $C$ be the connected
component of $T_2 \setminus \{\alpha\}$ containing $\epsilon$.
Since $\epsilon$ is an end of $\Lambda_1$ and $\epsilon$ is
contained in $C$, there must be infinitely many vertices in
$V\Lambda_1$ that lie in $C$. Choose $\alpha^\prime \in V\Lambda_1
\setminus \{\alpha\}$ such that $\alpha^\prime \in C$, and fix any
vertex $\delta \in [\alpha, \epsilon)_{T_2} \cap [\alpha^\prime,
\epsilon)_{T_2}$. Let $\beta$ be the vertex adjacent to $x$ in
$[x, \delta]_{T_1}$. Since $\beta$ is adjacent to $x$ in $T_1$, we
have $\beta \in V\Lambda_1$. If $\beta = \alpha$ then redefine
$\alpha := \alpha^\prime$. In this way, we can be sure $\beta \not
= \alpha$.

Since $\beta \in V\Lambda_1$ and $G_{\{\Lambda_1\}}$ is primitive
but not regular on $V\Lambda_1$, there is an element $g \in
G_{\alpha, \{\Lambda_1\}}$ such that $\beta^g \not = \beta$. Now
$\beta \in [\alpha, \delta]_{T_1}$, so $[\alpha, \delta]_{T_1}
\not = [\alpha, \delta]_{T_1}^g$. However, $g \in G_\alpha$, so
$[\alpha, \delta]_{T_1}^g = [\alpha, \delta^g]_{T_1}$ and
therefore $\delta \not = \delta^g$. There is only one half-line in
$T_1$ with initial vertex $\alpha$ contained in the end
$\epsilon$, so if $g \in G_{\alpha, \epsilon}$ then $g$ must fix
the half-line $[\alpha, \epsilon)_{T_1}$ pointwise. The vertex
$\delta$ lies on this half-line, so $G_{\alpha, \epsilon} \leq
G_\delta$. Since $\delta \not = \delta^g$, the automorphism $g$
thus cannot lie in $G_\epsilon$. This, however, is absurd, as
$\epsilon$ is the end of $\Lambda_1$, and $g \in G_{\alpha,
\{\Lambda_1\}} \leq G_{\alpha, \epsilon}$. \end{proof}

If $\Gamma$ is a locally finite connected graph or digraph, an automorphism
$g \in \aut \Gamma$ is called {\it hyperbolic} if $g$ fixes
precisely two thin ends of $\Gamma$ and leaves no non-empty finite
subset of $V\Gamma$ invariant.

\begin{thm} {\normalfont(\cite[Proposition 3.4]{tits1970})} \label{thm:tits}
If $T$ is an undirected tree and $G$ is a group of automorphism of
$T$ containing no hyperbolic elements, then either $G$ fixes some
vertex or leaves some edge invariant, or $G$ fixes some end of the
tree $T$. \qed \end{thm}

\begin{lemma} \label{lemma:centroid}
Suppose $G$ is a primitive group of permutations of an infinite
set and $\Gamma_1$ and $\Gamma_2$ are canonical orbital digraphs of
$G$. Let $T_2$ be the block-cut-vertex tree of $\Gamma_2$. If
$\Lambda_1$ is a lobe of $\Gamma_1$ then there exists $x \in
VT_2$ such that $G_{\{\Lambda_1\}} \leq G_{x}$, and for all
$\alpha_1, \alpha_2 \in V\Lambda_1$ we have $d_{T_2}(\alpha_1, x)
= d_{T_2}(\alpha_2, x)$. \end{lemma}

\begin{proof}
If $\Lambda_1$ is infinite, then it contains exactly one lobe-end
$\epsilon$ of $\Gamma_1$, which by
Lemma~\ref{lemma:thick_ends_go_to_thick_ends} is a lobe-end of
$\Gamma_2$. Let $\Lambda_2$ be the lobe of $\Gamma_2$ whose end
is $\epsilon$ and let $x$ be the vertex of $T_2$ corresponding to
the lobe $\Lambda_2$. Then $G_{\{\Lambda_1\}} = G_{\epsilon} =
G_{\{\Lambda_2\}} = G_x$, so we have $G_{\{\Lambda_1\}}$ fixes
some vertex in $T_2$.

Now assume $V\Lambda_1$ is finite. Since $G_{\{\Lambda_1\}}$
cannot contain a hyperbolic element, we may use
Theorem~\ref{thm:tits} to deduce that $G_{\{\Lambda_1\}}$ must fix
a vertex of $T_2$, or leave some edge of $T_2$ invariant, or fix
an end of $T_2$. In the latter two cases, the group must still fix
some vertex of $T_2$. Indeed, since no element of $G$ may
interchange two adjacent vertices of $T_2$, if the group
$G_{\{\Lambda_1\}}$ leaves some edge $\{\alpha, x\}$ of $T_2$
invariant, then it must fix both $\alpha$ and $x$ pointwise.
Furthermore, if $\epsilon$ is some end of $T_2$ fixed by the group
$G_{\{\Lambda_1\}}$, take $x$ to be any vertex in $T_2$ such that
the connected component of $T_2 \setminus \{x\}$ containing the
end $\epsilon$ contains no element in $V\Lambda_1$. Since
$V\Lambda_1$ is finite, there are infinitely many choices for such
a vertex. Then any half-line in $\epsilon$ whose root lies in
$V\Lambda_1$ must contain the vertex $x$. Thus, if
$G_{\{\Lambda_1\}}$ fixes the end $\epsilon$, it must also fix the
vertex $x$.

Finally we observe that, since $G_{\{\Lambda_1\}}$ is transitive
on $V\Lambda_1$, if $\alpha_1, \alpha_2 \in V\Lambda_1$ then
$d_{T_2}(\alpha_1, x) = d_{T_2}(\alpha_2, x)$. \end{proof}

Henceforth, if $\Lambda_1$ is a lobe of $\Gamma_1$, the vertex
$x$ in the above lemma will be referred to as the {\it centroid of
$\Lambda_1$ in $T_2$}. Since the above arguments are symmetric in
$\Gamma_1$ and $\Gamma_2$, we can also make reference to the {\it
centroid of $\Lambda_2$ in $T_1$}, where $\Lambda_2$ is any lobe
of $\Gamma_2$, and $T_1$ is the block-cut-vertex tree of
$\Gamma_1$.

Let $\Omega$ be some fixed infinite set and let $G$ be a primitive
group of permutations of $\Omega$. Given any locally
finite primitive connectivity-one orbital digraph $\Gamma_1$ of $G$,
if each vertex lies in $m$ distinct lobes and $\Lambda$ is a
lobe of $\Gamma_1$, we write the digraph $\Gamma_1$ as $\Gamma(m,
\Lambda)$. By Theorem~\ref{thm:connectivity1}, all locally finite primitive connectivity-one orbital
digraphs of $G$ can be written in this way. Since $G$ acts
transitively on the set of lobes of $\Gamma(m, \Lambda)$, there
is a natural equivalence relation on the set of such digraphs: the
orbital digraph $\Gamma_2$ is equivalent to $\Gamma_1$ if $\Gamma_2
= \Gamma(m, \Lambda^\prime)$ for some lobe $\Lambda^\prime$ of
$\Gamma_2$ satisfying $V\Lambda = V\Lambda^\prime$. It is under
this equivalence relation that a canonical orbital digraph of a
primitive permutation group can be considered essentially unique.

\begin{thm} \label{thm:structure_tree_is_unique}
If $G$ is a primitive group of permutations of an infinite set
$\Omega$ possessing no infinite subdegree, then the canonical
orbital digraphs of $G$ are equivalent.
\end{thm}

\begin{proof}
The group $G$ has at least one canonical orbital digraph by
Theorem~\ref{thm:connectivity1}. Suppose $G$ has two canonical
orbital digraphs $\Gamma_1 = \Gamma(m_1, \Lambda_1)$ and $\Gamma_2 =
\Gamma(m_2, \Lambda_2)$. By
Theorem~\ref{thm:directed_connectivity_one_graphs}, both $m_1$ and
$m_2$ are at least two, and the digraphs $\Lambda_1$ and $\Lambda_2$
are primitive but not automorphism-regular, and have at least
three vertices. Let $T_1$ and $T_2$ be the block-cut-vertex trees
of $\Gamma_1$ and $\Gamma_2$ respectively. We wish to show there
is a lobe $\Lambda_2^\prime$ of $\Gamma_2$ such that $V\Lambda_1
= V\Lambda_2^\prime$. From this we will deduce that the digraphs
$\Gamma_1$ and $\Gamma_2$ are equivalent. Set $d:= \mathrm{sup}
\{d_{T_2}(\alpha, \beta) \mid \alpha, \beta \in V\Lambda_1\}$.
This is finite, by Lemma~\ref{lemma:centroid}, so we may fix
$\alpha, \beta \in V\Lambda_1$ such that $d_{T_2}(\alpha, \beta) =
d$. Let $x_2$ be the vertex of $T_2$ that is adjacent to $\alpha$
in the line $[\alpha, \beta]_{T_2}$. Since $x_2$ is adjacent to
$\alpha$ in $T_2$, it must lie in $VT_2 \setminus V\Gamma_2$, and
therefore corresponds to some lobe $\Lambda_2^\prime$ of
$\Gamma_2$. Let $c_1$ be the centroid of $\Lambda_1$ in $T_2$.

We begin by showing $G_{\alpha, \{\Lambda_1\}} = G_{\alpha,
\{\Lambda_2^{\prime}\}}$. Now $c_1$ and $\beta$ must lie in the
same connected component of $T_2 \setminus \{\alpha\}$, and
$\beta$ and $x_2$ lie in the same component of $T_2 \setminus
\{\alpha\}$, so $c_1$ and $x_2$ are in the same component of $T_2
\setminus \{\alpha\}$. Since $x_2$ is adjacent to $\alpha$ in
$[\alpha, \beta]_{T_2}$, it follows that $x_2 \in [\alpha,
c_1]_{T_2}$, and therefore $G_{\alpha, c_1} \leq G_{x_2}$. Thus
\[G_{\alpha, \{\Lambda_1\}} \leq G_{\alpha, c_1} \leq G_{\alpha, x_2} = G_{\alpha, \{\Lambda_2^\prime\}}.\]

Now choose $\gamma \in \Lambda_2^\prime \setminus \{\alpha\}$ and
repeat this argument, replacing the vertex $\beta$ with $\gamma$,
the tree $T_2$ with $T_1$ and the lobe $\Lambda_1$ with
$\Lambda_2^\prime$. Whence, there exists a lobe
$\Lambda_1^\prime$ of $\Gamma_1$ containing $\alpha$ such that
$G_{\alpha, \{\Lambda_2^\prime\}} \leq G_{\alpha,
\{\Lambda_1^\prime\}}$, and therefore
\[G_{\alpha, \{\Lambda_1\}} \leq G_{\alpha, \{\Lambda_2^\prime\}} \leq G_{\alpha, \{\Lambda_1^\prime\}}.\]
Since $G_\alpha$ is transitive on the lobes of $\Gamma_1$
containing $\alpha$, there exists $g \in G_\alpha$ such that
$\Lambda_1^g = \Lambda_1^{\prime}$. As there are only finitely
many such lobes, $\Lambda_1^{g^n} = \Lambda_1$ for some natural
number $n$. Now $G_{\alpha, \{\Lambda_1\}} \leq G_{\alpha,
\{\Lambda_1\}^g}$, and therefore ${G_{\alpha, \{\Lambda_1\}^g}}
\leq G_{\alpha, \{\Lambda_1\}^{g^2}}$. If we continue in this way
we eventually obtain
\[G_{\alpha, \{\Lambda_1\}} \leq G_{\alpha, \{\Lambda_1\}^g} \leq \cdots \leq G_{\alpha, \{\Lambda_1\}^{g^n}};\]
however, $\Lambda_1^{g^n} = \Lambda_1$, so it follows that
$G_{\alpha, \{\Lambda_1\}} = G_{\alpha, \{\Lambda_1\}^g}$. Thus
\[G_{\alpha, \{\Lambda_1\}} = G_{\alpha, \{\Lambda_2^{\prime}\}}.\]

We now show $G_{\{\Lambda_1\}} = G_{\{\Lambda_2^{\prime}\}}$. Let
$\gamma$ be the vertex in $[x_2, \beta]_{T_2}$ adjacent to $x_2$,
and observe $\gamma \in V\Lambda_2^\prime$. Now
$G_{\{\Lambda_2^\prime\}}$ does not act regularly on
$V\Lambda_2^\prime$, so there is an element $g \in G_{\alpha,
\{\Lambda_2^\prime\}}$ such that $\gamma^g \not = \gamma$.
Therefore $x_2$ is the only vertex in $[x_2, \beta]_{T_2} \cap
[x_2, \beta^g]_{T_2}$. Thus $d_{T_2}(\beta, \beta^g) =
d_{T_2}(\beta, x_2) + d_{T_2}(\beta^g, x_2) = (d-1) + (d-1) =
2d-2$. Now $g \in G_{\alpha, \{\Lambda_2^\prime\}} = G_{\alpha,
\{\Lambda_1\}}$, and $\beta \in V\Lambda_1$, so $\beta^g \in
V\Lambda_1$. Furthermore, $d = d_{T_2}(\alpha, \beta)$ must be
even, so if $d > 2$, then $d \geq 4$, and therefore
$d_{T_2}(\beta, \beta^g) = 2d-2 > d$. This is not possible, as $d$
was chosen to be maximal. It must therefore be the case that
$d=2$.

Since no two vertices of $\Lambda_1$ are at distance
greater than $d$ in $T_2$, all vertices of $\Lambda_1$ are
adjacent to $x_2$ in $T_2$, and therefore lie in
$\Lambda_2^\prime$. Now $\Lambda_1$ contains at least three
vertices, each of which is adjacent to $x_2$ in $T_2$, and so,
since $G_{\{\Lambda_1\}}$ is transitive on $V\Lambda_1$, it must
be the case that $G_{\{\Lambda_1\}}$ fixes the vertex $x_2$, and
therefore fixes $V\Lambda_2^\prime$ setwise. Thus
$G_{\{\Lambda_1\}} \leq G_{\{\Lambda_2^\prime\}}$. If
$G_{\{\Lambda_1\}} \not = G_{\{\Lambda_2^\prime\}}$ then
$G_{\alpha, \{\Lambda_2^\prime\}} = G_{\alpha, \{\Lambda_1\}} <
G_{\{\Lambda_1\}} < G_{\{\Lambda_2^\prime\}}$, which cannot happen
since $G_{\{\Lambda_2^\prime\}}$ is primitive on $V\Lambda_2$.
Hence $G_{\{\Lambda_1\}} = G_{\{\Lambda_2^\prime\}}$.

Finally, since $G_{\{\Lambda_1\}}$ and $G_{\{\Lambda_2^\prime\}}$
are transitive on $V\Lambda_1$ and $V\Lambda_2^\prime$
respectively, and $\alpha \in V\Lambda_1 \cap V\Lambda_2^\prime$,
we have $V\Lambda_2^\prime = \alpha^{G_{\{\Lambda_2^\prime\}}} =
\alpha^{G_{\{\Lambda_1\}}} = V\Lambda_1$.

Now, $\Lambda_1$ is vertex-transitive, so the out-valency and
in-valency of each vertex must be non-zero. Since $\Gamma_1$
is arc-transitive, this implies that $G_\alpha$
is transitive on the set of lobes of $\Gamma_1$ that
contain $\alpha$. Similarly, $G_\alpha$ is transitive
on the set of lobes of $\Gamma_2$
containing $\alpha$. Whence, $m_1 = |G_\alpha:G_{\alpha,
\{\Lambda_1\}}| = |G_{\alpha}:G_{\alpha, \{\Lambda_2^\prime\}}| =
m_2$. Therefore $\Gamma_2 = \Gamma(m_1, \Lambda_2^\prime)$ and
$V\Lambda_2^\prime = V\Lambda_1$, so the digraphs $\Gamma_1$ and
$\Gamma_2$ are equivalent. \end{proof}

It is now possible to describe every locally finite
connectivity-one orbital digraph of a primitive group $G$.

\begin{thm} \label{thm:list_of_k1_graphs}
If $G$ is a primitive group of permutations of an infinite set
$\Omega$ with more than one permutation-end and no infinite
subdegree, then $G$ has a canonical orbital digraph $\Gamma(m,
\Lambda)$, and the canonical orbital digraphs of $G$ are precisely
the digraphs $\Gamma(m, \Lambda^\prime)$, where $\Lambda^\prime$ is
an orbital digraph of $G_{\{\Lambda\}}$ acting on the set
$V\Lambda$. \end{thm}

\begin{proof}
Suppose $G$ has a locally finite orbital digraph with more than one
end. By Theorem~\ref{thm:connectivity1}, $G$ has a canonical
orbital digraph $\Gamma(m, \Lambda)$. Set $\Gamma:= \Gamma(m,
\Lambda)$, and fix vertices $\alpha, \beta \in V\Lambda$ such that
$\Gamma = (\Omega, (\alpha, \beta)^G)$. Suppose $\Gamma^\prime$ is
also a canonical orbital digraph of $G$. We show there exists an
orbital digraph $\Lambda^\prime$ of $G_{\{\Lambda\}}$ acting on
$V\Lambda$ such that $\Lambda^\prime$ is a lobe of
$\Gamma^\prime$ and $\Gamma^\prime = \Gamma(m, \Lambda^\prime)$.
By Theorem~\ref{thm:structure_tree_is_unique}, the digraph
$\Gamma^\prime$ must be equivalent to $\Gamma$, so there exists a
lobe $\Lambda^\prime$ of $\Gamma^\prime$ such that $\Gamma^\prime
= \Gamma(m, \Lambda^\prime)$ and $V\Lambda^\prime = V\Lambda$. Let
$(\alpha, \beta^\prime)$ be an arc in $\Lambda^\prime$. Since
$\Gamma^\prime$ is arc-transitive, $\Gamma^\prime = (\Omega,
(\alpha, \beta^\prime)^G)$. If $(\gamma, \delta)$ is any arc in
$\Lambda^\prime$, then there exists an element $g \in G$ such that
$(\gamma, \delta) = (\alpha, \beta^\prime)^g$, but such an
automorphism must fix the lobe $\Lambda^\prime$, and therefore
lies in $G_{\{\Lambda^\prime\}} = G_{\{\Lambda\}}$. Hence
$\Lambda^\prime$ is an orbital digraph of $G_{\{\Lambda\}}$ acting
on the set $V\Lambda$.

Conversely, suppose $\Lambda^\prime$ is an orbital digraph of
$G_{\{\Lambda\}}$ acting on $V\Lambda$. We show the
connectivity-one digraph $\Gamma(m, \Lambda^\prime)$ is a canonical
orbital digraph of $G$.

It is simple to check that $G \leq \aut \Gamma(m, \Lambda')$, so
we begin by showing $G$ is arc-transitive on $\Gamma(m,
\Lambda')$. Suppose $\Lambda'_1$ and $\Lambda'_2$ are lobes of
$\Gamma(m, \Lambda^\prime)$. Then there exist lobes $\Lambda_1$
and $\Lambda_2$ of $\Gamma(m, \Lambda)$ such that $V\Lambda_1 =
V\Lambda'_1$ and $V\Lambda_2 = V\Lambda'_2$. Since $G$ acts
arc-transitively on $\Gamma(m, \Lambda)$, it acts transitively on
its lobes. Therefore, there exists an automorphism $g \in G$ such
that $\Lambda_1^g = \Lambda_2^g$. Since $G \leq \aut \Gamma(m,
\Lambda')$, we must also have $\Lambda^{\prime \ g}_1 =
\Lambda'_2$, so $G$ permutes the lobes of $\Gamma(m, \Lambda')$
transitively. Thus $G$ is transitive on the lobes of $\Gamma(m,
\Lambda)$ and acts arc-transitively on each lobe, so it must act
arc-transitively on the whole digraph $\Gamma(m, \Lambda')$.

It remains to show the lobes of $\Gamma(m, \Lambda')$ have at
most one end. This follows immediately from
Lemma~\ref{lemma:orbital_graphs_have_same_ends}, since both
$\Lambda$ and $\Lambda'$ are orbital digraphs of $G_{\{\Lambda\}}$
acting on $V\Lambda$, they must have the same ends.

Hence, $\Gamma(m, \Lambda)$ is a connectivity-one orbital digraph of
$G$ in which each lobe has at most one end, and is therefore a
canonical orbital digraph of the group $G$. \end{proof}

Thus, for any primitive group $G$ possessing a canonical orbital
digraph $\Gamma$, when speaking of the permutation-ends of $G$ as
being the ends of $\Gamma$, it makes sense to refer to the
{\em lobe-ends} and {\em tree-ends of $G$} as being, respectively,
the lobe-ends and tree-ends of its canonical orbital digraph. This
can be extended to include any group $G$ with just one permutation-end
$\epsilon$, by defining the {\em lobe-end of $G$} to be the end $\epsilon$.\\

That a whole class of infinite primitive permutation groups
possessing no infinite subdegree should each have an orbital digraph
whose structure is tree-like is a property that will be explored
and exploited in a future paper by the author. In particular, this
result will be fundamental in the examination of the possible
rates of subdegree growth of infinite primitive permutation
groups.\\

The concluding section of this paper is conceived as a simple
illustration of just one way in which the existence of the
canonical orbital digraph may be used to easily obtain a deep
understanding of the structure of infinite primitive groups
possessing no infinite subdegree, through the application of
Bass and Serre's theory on trees.

\section{The structure of infinite primitive groups}
\label{section:structure_of_G}

A common example of an amalgamated free product is when
there are two groups $G_1$ and $G_2$, with subgroups $H_1$ and
$H_2$ that are isomorphic via $\varphi : H_1 \rightarrow H_2$. In
this case the free product of $G_1$ and $G_2$ with amalgamated
subgroup $H_1$ is denoted by $G_1 \ast_{H_1} G_2$.

Let $G$ be a group acting on a digraph $\Gamma$. An inversion is a
pair consisting of an element $g \in G$ and an arc $(x,
y) \in A\Gamma$ such that $(x, y)^g = (y, x)$; if there is no such
pair the group $G$ is said to act {\em without inversion} on
$\Gamma$. If $G$ acts on $\Gamma$ without inversion, we define the
{\em quotient digraph} $G / \Gamma$ to be the digraph whose vertex and
arc set are the quotients of $V\Gamma$ and $A\Gamma$ respectively
under the action of $G$; a {\em fundamental domain} of $\Gamma
\gmod G$ is a subgraph of $\Gamma$ that is isomorphic to $G /
\Gamma$. A {\em segment} of $\Gamma$ is a subgraph of $\Gamma$
consisting two adjacent vertices and an arc of $\Gamma$
between them.

\begin{thm} {\normalfont(\cite[Theorem 6]{serre:trees})} \label{thm:serre:G_on_T_is_free} Let $G$ be a group acting on a digraph $\Gamma$ and let $P = \{\{x, y\}, (x, y)\}$ be a segment of $\Gamma$. If $P$ is a fundamental domain of $\Gamma \gmod G$ then $\Gamma$ is a tree if and only if $G \cong G_x \ast_{G_{x, y}} G_y$. \qed \end{thm}

Using this powerful theorem, we may determine the structure of
infinite primitive permutation groups with more
than one permutation-end and no infinite subdegree.

\begin{thm} \label{thm:G_amalgamated_free_product} If $G$ is an infinite primitive group of permutations of an infinite set $\Omega$ with more than one
permutation-end and no infinite subdegree, then $G$ has a canonical orbital digraph
$\Gamma(m, \Lambda)$ and
\[G \cong G_{\alpha} \ast_{G_{\alpha, \{\Lambda\}}} G_{\{\Lambda\}},\]
where $\alpha \in V\Lambda$, and $G_{\alpha, \{\Lambda\}}$ is a
non-trivial maximal proper subgroup of $G_{\{\Lambda\}}$ that
fixes no element in $V\Lambda \setminus \{\alpha\}$. \end{thm}

\begin{proof} Let $G$ be an infinite primitive group possessing a finite suborbit whose pair is also finite, and suppose $G$ has an orbital digraph with more than one end. By Theorem~\ref{thm:connectivity1}, $G$ has a canonical orbital digraph $\Gamma$ of the form $\Gamma(m, \Lambda)$, where $m \geq 2$ and $\Lambda$ is a primitive but not automorphism-regular digraph, with at least three vertices and at most one end. Let $T$ be the block-cut-vertex tree of $\Gamma$. Then $G \leq \aut T$. There is a natural bipartition of $T$, one part containing the vertices of $\Gamma$, and the other containing vertices corresponding to the lobes of $\Gamma$. Since $G$ preserves this bipartition, the group $G$ acts on $T$ without inversion, and has two orbits on the vertices of $T$.

Let $x \in VT$ correspond with the lobe $\Lambda$ of $\Gamma$,
and fix $\alpha \in V\Lambda$. Then $(\alpha, x)$ is an arc in
$T$. If we define $T^\prime := (VT, (\alpha, x)^G)$, then
$T^\prime$ is an arc-transitive tree upon which $G$ acts without
inversion, with fundamental domain the segment $\{\{\alpha, x\},
(\alpha, x)\}$. Hence, by Theorem~\ref{thm:serre:G_on_T_is_free},
\[G \cong G_{\alpha} \ast_{G_{\alpha, x}} G_x.\]
Since $\Lambda$ is the lobe of $\Gamma$ corresponding to the
vertex $x$ of $T$, the stabiliser $G_x$ is equal to
$G_{\{\Lambda\}}$. As $G_{\{\Lambda\}}$ acts primitively on
$V\Lambda$, the group $G_{\alpha, \{\Lambda\}}$ is maximal in
$G_{\{\Lambda\}}$. By
Theorem~\ref{theorem:block_stabiliser_not_regular}, $G_{\alpha,
\{\Lambda\}}$ fixes no vertex in $V\Lambda \setminus \{\alpha\}$.
\end{proof}


This paper forms part of the author's DPhil thesis, completed
at the University of Oxford, under the supervision of Peter
Neumann. The author would like to thank Dr Neumann for his help and guidance.
The author would also like to
thank the EPSRC for generously funding this research.


\end{document}